\documentclass[a4paper]{article}

\usepackage{url}
\usepackage{amsmath}
\usepackage{amsfonts}
\usepackage{graphicx}
\usepackage{epstopdf}
\usepackage{algorithmic}
\ifpdf
  \DeclareGraphicsExtensions{.eps,.pdf,.png,.jpg}
\else
  \DeclareGraphicsExtensions{.eps}
\fi

\usepackage{cite}
\usepackage{amssymb}
\newcommand{\J}{\mathrm{j}}  
\newcommand{\D}{\mathrm{d}}  
\newcommand{\E}{\mathrm{e}}  
\newcommand{\rbar}{\overline{r}}
\newcommand{\rmid}{r_{\text{mid}}}
\newcommand{\rmax}{r_{\text{max}}}
\newcommand{\Rmid}{R_{\text{mid}}}
\newcommand{\Rbar}{\overline{R}}

\title{Closed-form evaluation of potential integrals in the Boundary
  Element Method \footnote{This work was partially carried out under
    the Horizon~2020 project AERIALIST, project identifier~723367}}

\author{Michael Carley\thanks{Department of Mechanical Engineering,
    University of Bath, Bath, BA2 7AY, United Kingdom}
  (\texttt{m.j.carley@bath.ac.uk})}

\usepackage{amsopn}

\begin{document}

\maketitle

\begin{abstract}
  A method is presented for the analytical evaluation of the singular
  and near-singular integrals arising in the Boundary Element Method
  solution of the Helmholtz equation. An error analysis is presented
  for the numerical evaluation of such integrals on a plane element,
  and used to develop a criterion for the selection of quadrature
  rules. The analytical approach is based on an optimized expansion of
  the Green's function for the problem, selected to limit the error to
  some required tolerance. Results are presented showing accuracy to
  tolerances comparable to machine precision. 
\end{abstract}

\section{Introduction}
\label{sec:intro}

A central part of the Boundary Element Method (BEM) is the evaluation
of potential integrals, to compute the contribution of an element to
the potential field, or to the entries of the solution matrix. It is
thus a key factor in the accuracy and efficiency of any
implementation, and one which has attracted great interest over many
decades. In this paper we develop a method for the evaluation of
integrals which arise in the three-dimensional BEM for the wave
equation, in particular in acoustics, where the acoustic potential
$\phi$, external to a surface $A$, is given by the integral
formulation:
\begin{align}
  \label{equ:potential}
  \phi(\mathbf{x}) &= 
  \int_{A} 
  \frac{\partial\phi_{1}}{\partial n}G(\mathbf{x},\mathbf{x}_{1})
  -
  \frac{\partial G(\mathbf{x},\mathbf{x}_{1})}{\partial n}\phi_{1}
  \,\D A,
\end{align}
where $\mathbf{x}$ indicates position, subscript $1$ variables of
integration on the surface $A$, and $n$ the outward pointing normal to
the surface. The Green's function $G$ is:
\begin{align}
  \label{equ:greens}
  G(\mathbf{x};\,\mathbf{x}_{1}) &= \frac{\E^{\J k R}}{4\pi R},\\
  R &= |\mathbf{x}-\mathbf{x}_{1}|, \nonumber
\end{align}
where $k$ is acoustic wavenumber. Given the surface potential $\phi$
and gradient $\partial\phi/\partial n$, the potential, and, after
differentiation, its gradient(s), can be evaluated at any point in the
field. Also, given a boundary condition for $\phi$ and/or
$\partial\phi/\partial n$ on $A$, the integral equation can be solved
for $\phi(\mathbf{x})$ and/or $\partial\phi/\partial n(\mathbf{x})$,
$\mathbf{x}\in A$.

If the boundary integral equation is solved using a collocation
method, the surface $A$ is divided into elements, here taken to be
plane triangles, and suitable shape functions are used to interpolate
the potential on these elements. The integral equation is transformed
to a linear system in the element potentials, with the influence
coefficients determined by the potential generated by each element at
each node of the surface mesh. This leads to the requirement to
evaluate integrals $I$ and $\partial I/\partial n$ where:
\begin{align}
  \label{equ:element}
  I &= \iint_{A_{e}} f(\xi,\eta)
  G(\mathbf{x},\mathbf{x}_{1}(\xi,\eta))\,\D A_{e},
\end{align}
with $A_{e}$ the surface of an element and $(\xi,\eta)$ a coordinate
system local to $A_{e}$. The requirement then is to evaluate integrals
of $\exp[\J k R]/R$ and its derivatives over a triangular
element. This is especially challenging when the field point
$\mathbf{x}$ is on, or near, the element, and the $1/R$ singularity
must be accommodated in the integration scheme. 

There are numerous numerical schemes for the evaluation of the surface
integrals, which mainly vary in their approach to dealing with the
singularity. There are also a number of analytical schemes for the
equivalent integral in the Laplace equation~\cite[for
example]{carley13,newman86,mogilevskaya-nikolskiy14,%
  carini-salvadori02,salvadori10,suh00}, some of which can be used to
deal with the singular terms in the acoustic problem and thus ease
numerical integration, but there are few analytical methods for the
Helmholtz problem. Clearly, given the absence of an exact analytical
solution for the retarded potential from a plane element, any
closed-form solution is an approximation, but it should be possible to
approximate the integral to any required accuracy, in a form amenable
to analytical manipulation, so that the result can be used as if it
were an analytical formula for the potential. This is especially
important for the case of a field point on or near the element, where
the ability to handle singularities analytically offers an advantage
over purely numerical schemes.

To the author's knowledge there are two published methods for
closed-form or analytical evaluation of the Helmholtz potential from a
planar
element~\cite{tadeu-antonio12,pourahmadian-mogilevskaya15}. These use
two different approaches to the problem. In
one~\cite{tadeu-antonio12}, an expression is derived in the Fourier
domain resulting in an expression based on a series of terms defined
by integrals of Hankel functions. These integrals can be evaluated
analytically in terms of Struve functions, yielding a closed-form
solution for the potential from a planar element, but at the expense
of using special functions not routinely available in numerical
libraries.

The second approach~\cite{pourahmadian-mogilevskaya15}, which is
similar in spirit to the method of this paper, makes use of results
derived for the Laplace problem~\cite{mogilevskaya-nikolskiy14} and
approximates $\exp[\J k R]$ as a polynomial over the element. This is
justified by noting that the element size is limited by the
requirement to avoid aliasing in the representation of the surface
potential, so that a relatively low-order approximation containing
five or six terms of the Taylor series for $\exp[\J k R]$ is adequate
for evaluation of the integrals to the tolerance specified.

The method of this paper uses a similar approach, in that it replaces
the exponential with an approximation of controlled error, and uses
results from an analysis of the Laplace problem~\cite{carley13} to
compute the terms in the resulting expansion. It differs in the form
of Laplace solution used, and in the choice of expansion for the
exponential, to give a systematic control on quadrature accuracy
optimized to require a minimum number of terms. Additionally, a
criterion is provided for choosing when, and when not, to use the
analytical approach or a purely numerical method, based on an error
analysis of integration using a polar coordinate transformation. To
the author's knowledge, this error analysis is novel and may have
applications more generally.

\section{Integration of the $1/R$ potential}
\label{sec:analysis:switch}

In order to motivate the development of the closed-form expression for
the acoustic potential, we begin by analyzing the numerical evaluation
of the Laplace potential, which corresponds to the leading-order,
singular, part of the Helmholtz potential, which gives rise to the
difficulties in numerical integration. 

\begin{figure}
  \centering
  \includegraphics{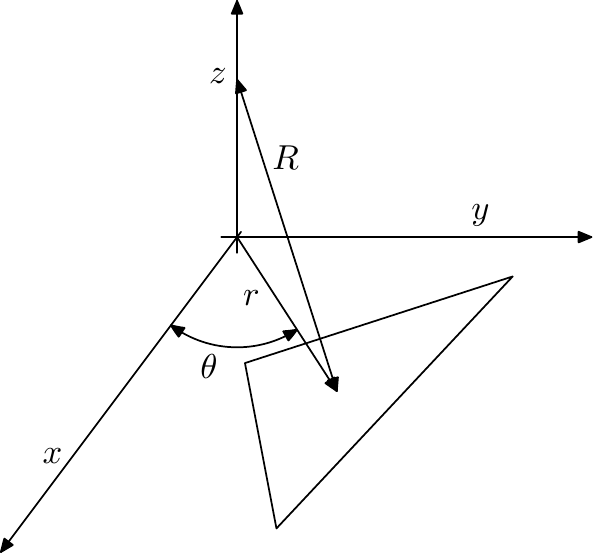}
  \caption{Integration of Laplace potential over a triangle}
  \label{fig:analysis:laplace}
\end{figure}

The model problem is shown in Figure~\ref{fig:analysis:laplace} and
consists of the evaluation of
\begin{align*}
  I &= \int_{A} \frac{r}{R}\,\D r\D\theta,\\
  R^{2} &= r^{2} + z^{2},
\end{align*}
over the area of the triangle shown, which lies in the plane $z=0$,
with the usual transformation to polar coordinates $(r,\theta)$ for
the integration.

The error in the evaluation of this integral, especially at small
values of $z$ arises from the singular, or near-singular, term
$1/R$. Here we develop an approximate error analysis for the
evaluation of this term, which can be used in determining the required
order of integration for $r/R$ or when to switch to some other
quadrature approach, such as that in the next section. An error
analysis for integration using the polar coordinate transformation has
been published previously~\cite{schwab-wendland92} but the analysis
presented here appears to be novel and is simple enough for use as an
\textit{a priori} estimator in determining quadrature order in
applications. 

The analysis depends on an error estimate for the $1/R$ term in a
numerical polar integration, such as~(\ref{equ:analysis:i0}). If such
an integration is performed using Gaussian quadrature, the integrand
is being approximated by a polynomial over the interval of integration
and the accuracy of the approximation is determined by the number of
terms required to approximate the integrand accurately. We perform the
analysis by estimating the error in the polynomial expansion of $1/R$
and use it to give an approximation of the order of polynomial
required to approximate $1/R$ to a given tolerance. Given this
polynomial order, a Gaussian quadrature of sufficiently high order can
be selected, or if the order required is too great, the analytical
method of the following section can be used.

\begin{figure}
  \centering
  \includegraphics{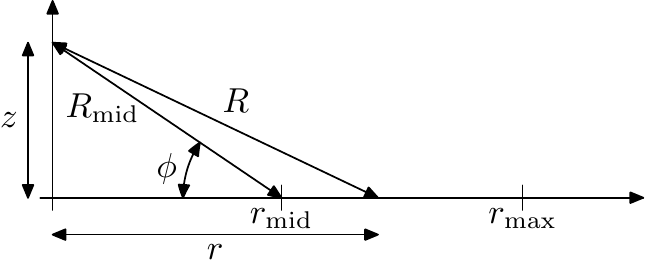}
  \caption{Notation for error analysis of $1/R$ expansion}
  \label{fig:error}
\end{figure}

From Figure~\ref{fig:error}, we write
\begin{align*}
  r &= \rmid - t\rmid,\quad \rmid=\rmax/2,\quad -1\leq t\leq 1,\\
  R^{2} &= r^{2} + z^{2} = \Rmid^{2}
  \left[
    1 - 2\cos^{2}\phi t + (t\cos\phi)^{2}
  \right],\\
  \Rmid^{2} &= \rmid^{2}+z^{2},\quad  \cos\phi = \rmid/\Rmid,
\end{align*}
and neglect the case of $\phi=0$ as in this case $r/R\equiv1$ and
the polynomial representation of the integrand raises no
difficulties.

Using the generating function for Legendre
polynomials~\cite[8.921]{gradshteyn-ryzhik80},
\begin{align}
  \frac{1}{R} &= 
  \frac{1}{\Rmid}
  \sum_{q=0}^{\infty}
  (t\cos\phi)^{q}P_{q}(\cos\phi).
\end{align}
If the expansion is truncated at $q=Q$, the error at any value of
$t$ is given by the remainder
\begin{align*}
  \epsilon_{Q}
  &=
  \frac{1}{R} - \frac{1}{\Rmid}
  \sum_{q=0}^{Q}
  (t\cos\phi)^{q}P_{q}(\cos\phi)
  =
  \frac{1}{\Rmid}
  \sum_{q=Q+1}^{\infty}
  (t\cos\phi)^{q}P_{q}(\cos\phi),
\end{align*}
which can be rewritten using the large-order asymptotic form of the
Legendre polynomial~\cite[8.918]{gradshteyn-ryzhik80},
\begin{align*}
  P_{q}(\cos\phi) \sim
  \left(
    \frac{2}{\pi q \sin\phi}
  \right)^{1/2}
  \cos
  \left[
    (q+1/2)\phi - \pi/4
  \right],
\end{align*}
so that
\begin{align*}
  \epsilon_{Q}
  &\approx
  \frac{1}{\Rmid}\frac{1}{(\pi\sin\phi)^{1/2}}
  \sum_{q=Q+1}^{\infty}
  t^{q}\frac{\cos^{q}\phi}{q^{1/2}}
  \left[
    \cos(q+1/2)\phi
    +
    \sin(q+1/2)\phi
  \right].
\end{align*}
An upper bound for the sum can be found by replacing $q^{1/2}$ with
$(Q+1)^{1/2}$ and, upon rearrangement,
\begin{align}
  \epsilon_{Q}
  &\approx
  \Im
  \frac{1+j}{\Rmid}\frac{\E^{\J\phi/2}}{[\pi(Q+1)\sin\phi]^{1/2}}
  \sum_{q=Q+1}^{\infty}
  \left(
    t\cos\phi\E^{\J\phi}
  \right)^{q},\nonumber\\
  &= 
  \Im
  \frac{1+j}{\Rmid}\frac{\E^{\J\phi/2}}{[\pi(Q+1)\sin\phi]^{1/2}}
  \frac{\left(t\cos\phi\E^{\phi}\right)^{Q+1}}{1-t\cos\phi\E^{\J\phi}}. 
  \label{equ:error}
\end{align}
As will be seen, this is an accurate estimate of the remainder in the
polynomial expansion of $1/R$ but it is oscillatory as a function of
$q$, so we adopt the more convenient measure of the magnitude rather
than the imaginary part,
\begin{align}
  \label{equ:error:abs}
  E_{Q} &= \frac{1}{\Rmid}
  \left(
    \frac{2}{\pi\sin\phi}
  \right)^{1/2}
  \frac{|t|^{Q+1}}{(Q+1)^{1/2}}
  \frac{\cos^{Q+1}\phi}{\left[(1-t)^{2}\cos^{2}\phi+\sin^{2}\phi\right]^{1/2}}.
\end{align}

We note that~(\ref{equ:error}) could be integrated over $t$ to give an
estimate of the total error in the integral of $r/R$ but this gives an
unwieldy expression with little advantage in implementations. Instead
we adopt as error criterion the absolute value given
by~(\ref{equ:error:abs}) with the value of $t$ given by the nearest
point on the element. In particular, when the projection of the field
point lies on the element, i.e. when the triangle
in Figure~\ref{fig:analysis:laplace} encloses the origin, $t=1$ and the
error estimate for $1/R$ is
\begin{align}
  \label{equ:error:abs:1}
  E_{Q} &= \frac{1}{\Rmid}
  \left(
    \frac{2}{\pi\sin^{3}\phi}
  \right)^{1/2}
  \frac{\cos^{Q+1}\phi}{(Q+1)^{1/2}}.
\end{align}
From the form of the error estimate, the reason for the difficulty in
evaluating near-singular integrals is clear: near the element plane
where $\phi\to0$, approximation of the integrand by a polynomial,
implicit in the use of Gaussian quadratures, incurs a very large
error, even for quite high order quadratures with large $Q$.

\begin{figure}
  \centering
  \includegraphics{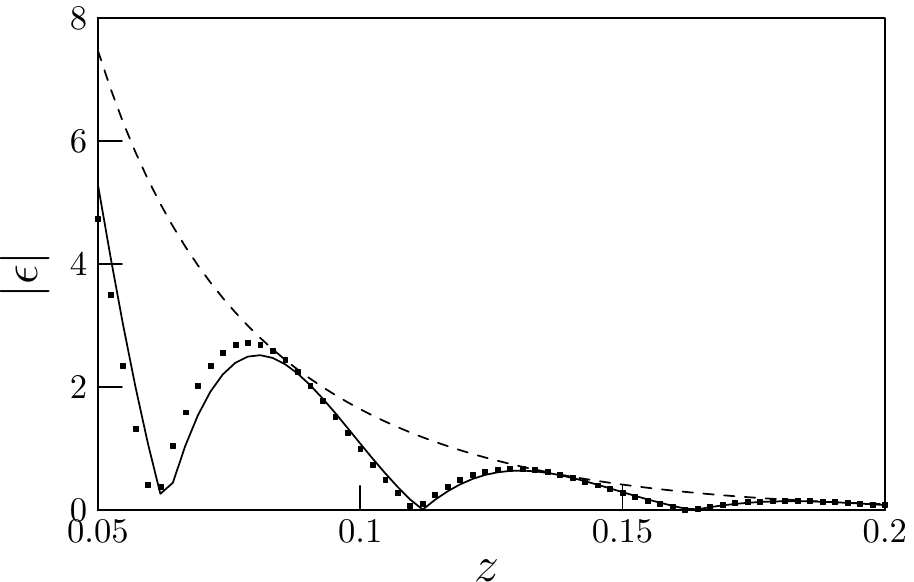}
  \caption{Error in polynomial approximation of $1/R$, $Q=32$,
    $\rmid=1/2$, $t=1$: solid line exact error; boxes estimate
    from~(\ref{equ:error}); dashed line absolute value
    from~(\ref{equ:error:abs})}
  \label{fig:error:check}
\end{figure}

\begin{figure}
  \centering
  \includegraphics{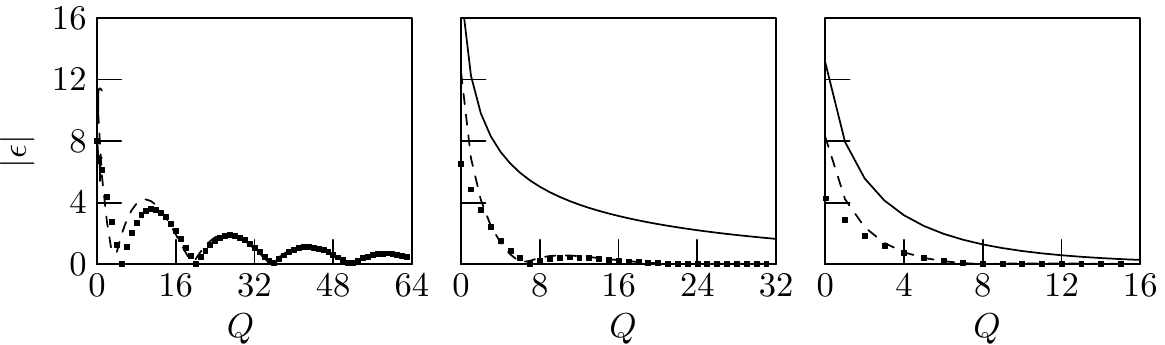}
  \caption{Error in polynomial approximation of $1/R$, $\rmid=1/2$,
    $z=0.1$, $t=1, 7/8, 3/4$: solid line exact error; boxes estimate
    from~(\ref{equ:error}); dashed line absolute value
    from~(\ref{equ:error:abs})}
  \label{fig:error:Q}
\end{figure}

To minimize the computational burden of using the criterion, it is
applied in the following manner. Given the transformation into
coordinates based on the element plane, the minimum distance
$r_{\min}$ from the origin to the triangle can be determined (see
Figure~\ref{fig:analysis:basic}) with $r_{\min}\equiv0$ when the
triangle encloses the origin. The maximum distance to a vertex
$r_{\max}$ is found similarly. Then we set $\rmid=r_{\max}/2$,
$t=(\rmid-r_{\min})/\rbar$ and other quantities as above. The
criterion is then applied by computing $E_{Q}$ for $Q=1,2,\ldots$
until $E_{Q}$ falls below some specified tolerance, and returning the
resulting value of $Q$, the order of polynomial required to compute
$1/R$ to the specified tolerance over the range of the integral. We
note that the error measure here is the maximum error in $1/R$ at any
point in the range of integration, which is quite a stringent, though
conservative, measure, but it will be found that $E_{Q}$ is a useful
assessment of the accuracy of quadrature.

Figure~\ref{fig:error:check} shows the error estimates as a function
of $z$ for a test case with a~32nd order polynomial, equivalent to
a~16 point Gaussian quadrature. The error estimate $\epsilon_{Q}$ is
seen to be very reliable, and the magnitude $E_{Q}$ does indeed match
the amplitude of the error quite closely. Figure~\ref{fig:error:Q}
shows the error as a function of $Q$ for fixed $z$ and again the error
behavior is accurately captured by the estimators. Despite the
relative simplicity of the error measures, they give reliable
indicators of the accuracy of the quadrature or of the order of
quadrature required for a given tolerance. We note finally that the
quantities used in the error measure are typically computed as part of
the geometric transformations required in generating a quadrature on
an element, so that there is very little overhead in applying the
error estimate.

\section{Analysis}
\label{sec:analysis}

\begin{figure}
  \centering
  \includegraphics{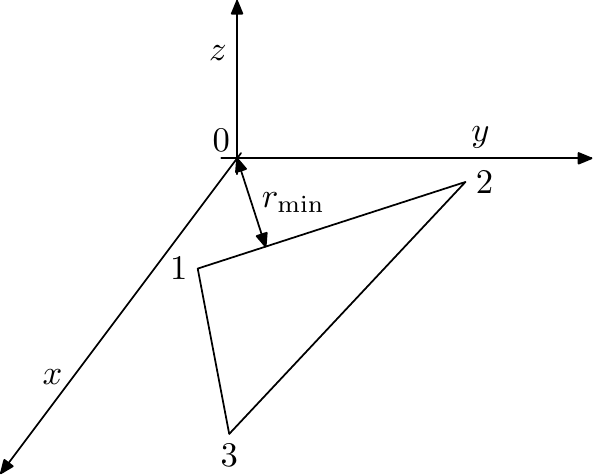}
  \includegraphics{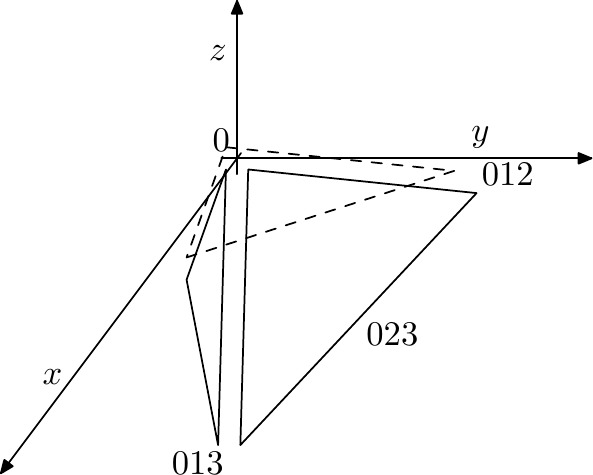} 
  \caption{Integration over a general triangle (left) by subdivision
    into three triangles centred at the origin (right). The triangle
    shown dashed in the exploded view on the right has negative
    orientation and its contribution is subtracted from that of the
    others. The distance $r_{\min}$ is used in applying the criterion
    of Section~\ref{sec:analysis:switch}.}
  \label{fig:analysis:basic}
\end{figure}

The problem to be considered is evaluation of the Helmholtz single-
and double-layer potential integrals on a planar triangular element.
Integration is performed after transformation of coordinates such that
the triangular element is defined by vertices $(x_{i},y_{i},0)$ and
the field point lies at $(0,0,z)$. The triangle is then decomposed
into up to three triangles each having a vertex at $(0,0,0)$. The
process is shown in Figure~\ref{fig:analysis:basic}. The approach is
similar to that taken in a previous analysis for the Laplace
potential~\cite{carley13}, though some changes are required to make it
suitable for the Helmholtz problem.

\begin{figure}[ht]
  \centering
  \includegraphics{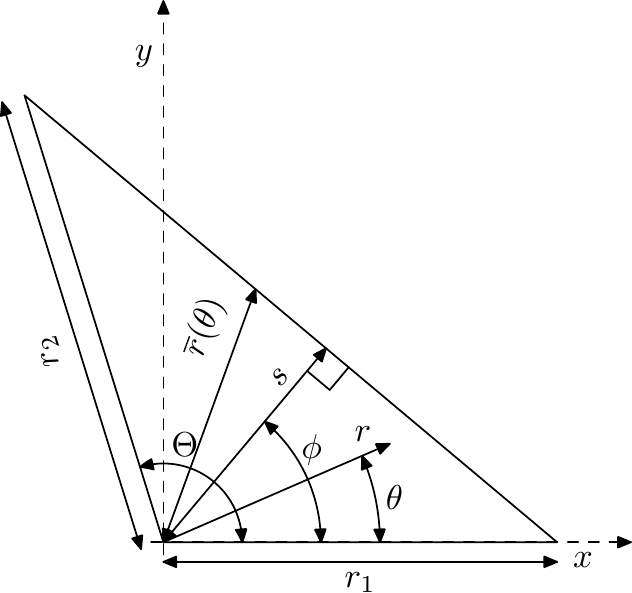}
  \caption{Reference triangle for integration}
  \label{fig:analysis:tri}
\end{figure}

Figure~\ref{fig:analysis:tri} shows the basic triangle which is used
for the evaluation of the contributions from the subtriangles of
Figure~\ref{fig:analysis:basic}. It has one vertex at the origin,
i.e.\ at the projection of the field point onto the element plane, and
is defined by the lengths of the two sides which meet at the origin,
$r_{1}$ and $r_{2}$, and by the angle $\Theta$ between them.

In developing the analysis, we assume that the triangular element
conforms to some reasonable standards of quality, in particular that
the edge length is no greater than some specified fraction of a
wavelength, typically between one sixth and one eighth. This
translates into a limit on $k\ell$, where $\ell$ is a typical edge
length. Taking into account the need to deal with triangles which are
larger than the element proper, such as triangle~$023$ in
Figure~\ref{fig:analysis:basic}, we assume that $k\ell<\pi/2$, which
allows us to limit the size of the expansions which will be employed
in evaluating the potential integrals, while retaining the required
accuracy. If necessary, this limit can be increased, at the expense of
extra computational effort and a small increase in stored data.

\subsection{Basic integrals}
\label{sec:analysis:basic}

Integration is performed on the reference triangle of
Figure~\ref{fig:analysis:tri}, using the polar coordinate system
$(r,\theta)$. Geometric parameters are defined,
\begin{align*}
  \phi &= \tan^{-1}
  \frac{r_{1}-r_{2}\cos\Theta}{r_{2}\sin\Theta},\,
  \rbar(\theta) = \frac{s}{\cos(\theta-\phi)},\\
  s &= r_{1}\cos\phi,\,
  S^{2} = s^{2}+z^{2},
\end{align*}
and auxiliary variables used in performing the integrations are
\begin{align*}
  \alpha^{2} = z^{2}/S^{2},\,
  \Rbar &= \left(\rbar^{2}+z^{2}\right)^{1/2} = S\Delta/\cos\theta,
  \quad
  \Delta^{2} = 1-\alpha^{2}\sin^{2}\theta.
\end{align*}

The integrals to be evaluated are the zeroth and first order
derivatives with respect to $z$ of
\begin{align*}
  I_{0} = \E^{\J k |z|}I_{0}',\quad 
  I_{x} = \E^{\J k |z|}I_{x}',\quad
  I_{y} = \E^{\J k |z|}I_{y}', 
\end{align*}
where the basic integrals are
\begin{subequations}
  \label{equ:analysis:i0xy}
  \begin{align}
    \label{equ:analysis:i0}
    I'_{0} &= 
    \int_{-\phi}^{\Theta-\phi}
    \int_{0}^{\rbar}
    \frac{\E^{\J k (R-|z|)}}{R}
    r\, \D r\,\D\theta,\\
    \label{equ:analysis:ix}
    I'_{x} &= 
    \int_{-\phi}^{\Theta-\phi}
    \int_{0}^{\rbar}
    \frac{\E^{\J k (R-|z|)}}{R}
    r^{2}\, \D r\,\cos\theta\,\D\theta,\\
    \label{equ:analysis:iy}
    I'_{y} &= 
    \int_{-\phi}^{\Theta-\phi}
    \int_{0}^{\rbar}
    \frac{\E^{\J k (R-|z|)}}{R}
    r^{2}\, \D r\,\sin\theta\,\D\theta,\\
    \text{with}\,R &= \left(r^{2} + z^{2}\right)^{1/2},\nonumber
  \end{align}  
\end{subequations}
and correspond to the zero and first order source terms required for
linear shape functions on a plane element. The normal derivatives are
given by differentiation with respect to $z$, for example,
\begin{align*}
  \frac{\partial I_{0}}{\partial z}
  &=
  \pm \J k \E^{\J k |z|}I_{0}' + \E^{\J k |z|}
  \frac{\partial I_{0}'}{\partial z}
  =
  -\frac{\partial I_{0}}{\partial n}
\end{align*}
where the element is oriented such that the normal lies in the
positive $z$ direction, and the upper (lower) signs are taken for
positive (negative) $z$. 

The integrals are evaluated by expanding the complex exponential in a
polynomial approximation, and evaluating term-by-term using analytical
formulae defined by recursion relations, as in previous
work~\cite{carley13}. The form of the approximation for $\exp[\J k x]$
will be considered in Section~\ref{sec:analysis:exp}, but for now we
write
\begin{align}
  \label{equ:exp:approx}
  \E^{\J k(R-|z|)} 
  &\approx 
  \sum_{q=0}^{Q}e_{q}
  k^{q}
  (R - |z|)^{q},\\
  e_{q} &= c_{q} + \J s_{q},\nonumber\\
  \text{where}\,\sin x &\approx \sum_{q=0}^{Q}s_{q}x^{q},\quad \cos x \approx
  \sum_{q=0}^{Q}c_{q}x^{q}.\nonumber
\end{align}
Expanding in powers of $k(R-|z|)$ has the advantages of ensuring that the
expansion remains valid for large values of $z$ as $(R-|z|)\to0$ as
$z\to\infty$, and providing a natural reduction in the number of terms
in~(\ref{equ:exp:approx}) for increasing $z$. 

Substituting (\ref{equ:exp:approx}) into (\ref{equ:analysis:i0xy}),
yields
\begin{subequations}
  \label{equ:analysis:ixy0:1}
  \begin{align}
  \label{equ:analysis:i0:1}
  I_{0}' &\approx 
  \sum_{q=0}^{Q} e_{q} K_{q,0},\\
  \label{equ:analysis:ix:1}
  I_{x}' &\approx 
  \sum_{q=0}^{Q} e_{q} K_{q,x},\\
  \label{equ:analysis:iy:1}
  I_{y}' &\approx 
  \sum_{q=0}^{Q} e_{q} K_{q,y},
\end{align}
\end{subequations}
where
\begin{subequations}
  \label{equ:analysis:k:1}
  \begin{align}
  \label{equ:analysis:k0:1}
    K_{q,0} &= k^{q}
    \int_{-\phi}^{\Theta-\phi}
    \int_{0}^{\rbar}
    \left(R - |z|\right)^{q}\frac{r}{R}
    \, \D r\,\D\theta,\\
    K_{q,x} &= k^{q}
    \int_{-\phi}^{\Theta-\phi}
    \int_{0}^{\rbar}
    \left(R - |z|\right)^{q}\frac{r^{2}}{R}
    \, \D r\cos\theta\,\D\theta,\\
    K_{q,y} &= k^{q}
    \int_{-\phi}^{\Theta-\phi}
    \int_{0}^{\rbar}
    \left(R - |z|\right)^{q}\frac{r^{2}}{R}
    \, \D r\sin\theta\,\D\theta.
  \end{align}
\end{subequations}
The integrals of (\ref{equ:analysis:k:1}) can be evaluated
analytically using a combination of recursions and tabulated
integrals. The inner integrals are given by,
\begin{align}
  k^{q}
  \int_{0}^{\rbar}
  \left(R - |z|\right)^{q}\frac{r}{R}
  \, \D r
  &=
  \frac{k^{q}(R - |z|)^{q+1}}{q+1}
  = 
  \frac{S(kS)^{q}}{q+1}
  \left(
    \frac{\Delta}{\cos\theta} - \alpha
  \right)^{q+1},\\
  k^{q}\int_{0}^{\rbar} \left(R - |z|\right)^{q}\frac{r^{2}}{R}\,
  \D r
  &=
  \frac{k^{q}\rbar}{q+2}\left(\Rbar - |z|\right)^{q+1}
  + 
  \frac{2|z|}{q+2} J_{q},\nonumber\\
  &= 
  \frac{sS(kS)^{q}}{q+2}
  \left(
    \frac{\Delta}{\cos\theta} - \alpha
  \right)^{q+1}
  \frac{1}{\cos\theta}
  + 
  \frac{2|z|}{q+2} J_{q},\\
  J_{q} &= 
   k^{q}\int_{(2|z|)^{1/2}}^{\left(\Rbar+|z|\right)^{1/2}}
  \left(t^{2}-2|z|\right)^{q+1/2}
  \,\D t.\nonumber
\end{align}
The integral $J_{q}$ can be evaluated using the recursion
\begin{align}
  J_{q}
  &= \frac{k^{q}\rbar}{2(q+1)}\left(\Rbar - |z|\right)^{q}
  - k|z| \frac{2q+1}{q+1}J_{q-1},\nonumber\\
  &= 
  \frac{s}{2}
  \frac{(kS)^{q}}{q+1}
  \left(
    \frac{\Delta}{\cos\theta} - \alpha
  \right)^{q}
  \frac{1}{\cos\theta}
  - k|z| \frac{2q+1}{q+1}J_{q-1},\\
  J_{0} &=
  \frac{\rbar}{2} -
  \frac{|z|}{2}\log\frac{\Rbar+\rbar}{|z|},\nonumber\\
  &=
  \frac{s}{2}\frac{1}{\cos\theta}
  +
  \frac{|z|}{4}\log\frac{\Delta-\alpha'}{\Delta+\alpha'},
\end{align}
so that all required terms are written in a form suitable for the
application of standard formulae for trigonometric
integrals~\cite[2.58]{gradshteyn-ryzhik80},
\begin{subequations}
  \begin{align}
    K_{q,0} &= 
    \frac{S(kS)^{q}}{q+1}
    \int_{-\phi}^{\Theta-\phi}
    \left(
      \frac{\Delta}{\cos\theta} - \alpha
    \right)^{q+1}
    \,\D\theta,\\    
    K_{q,x} &=
    \frac{sS(kS)^{q}}{q+2}
    \int_{-\phi}^{\Theta-\phi}
    \left(
      \frac{\Delta}{\cos\theta} - \alpha
    \right)^{q+1}
    \,\D\theta
    + 
    \frac{2|z|}{q+2} I_{q,c},\\
    K_{q,y} &=
    \frac{sS(kS)^{q}}{q+2}
    \int_{-\phi}^{\Theta-\phi}
    \left(
      \frac{\Delta}{\cos\theta} - \alpha
    \right)^{q+1}
    \frac{\sin\theta}{\cos\theta}
    \,\D\theta
    + 
    \frac{2|z|}{q+2} I_{q,s},\\
    I_{q,c} &= 
    \int_{-\phi}^{\Theta-\phi}
    J_{q}\cos\theta
    \,\D\theta,\,
    I_{q,s} = 
    \int_{-\phi}^{\Theta-\phi}
    J_{q}\sin\theta
    \,\D\theta.\nonumber
  \end{align}
\end{subequations}

The normal derivatives of the integrals can be evaluated by
differentiating terms, yielding,
\begin{align}
  \label{equ:analysis:dk:0}
  \frac{\partial K_{q,0}}{\partial z}
  &=
  \mp (kS)^{q}
  \int_{-\phi}^{\Theta-\phi}
  \left(
    \frac{\Delta}{\cos\theta} - \alpha
  \right)^{q+1}
  \frac{\cos\theta}{\Delta}
  \,\D \theta,
  \\
  \label{equ:analysis:dk:x}
  \frac{\partial K_{q,x}}{\partial z}
  &=
  \mp
  s
  (kS)^{q}\frac{q+1}{q+2}
  \int_{-\phi}^{\Theta-\phi}
  \left(
    \frac{\Delta}{\cos\theta} - \alpha
  \right)^{q+1}
  \frac{\cos\theta}{\Delta}
  \,\D\theta
  \pm \frac{2}{q+2}I_{q,c}
  + \frac{2|z|}{q+2}
  \frac{\partial I_{q,c}}{\partial z},\\
  \label{equ:terms:dz:1}
  \frac{\partial K_{q,y}}{\partial z}
  &=
  \mp
  s
  (kS)^{q}\frac{q+1}{q+2}
  \int_{-\phi}^{\Theta-\phi}
  \left(
    \frac{\Delta}{\cos\theta} - \alpha
  \right)^{q+1}
  \frac{\sin\theta}{\Delta}
  \,\D\theta
  \pm \frac{2}{q+2}I_{q,s}
  + \frac{2|z|}{q+2}
  \frac{\partial I_{q,s}}{\partial z}.
\end{align}
All integrals can then be evaluated using the results of
Section~\ref{sec:basic:integrals}. This gives a means of evaluating
all required expressions for the integrals on the triangular element,
which can then be summed to give the integral over the initial general
triangle.

\subsection{Hypersingular integral}
\label{sec:analysis:hyper}

The results of (\ref{sec:analysis:basic}) may be used to solve
boundary integral problems using a standard Helmholtz equation. It is
often desirable to employ a Burton and Miller
approach~\cite{burton-miller71} to avoid the well-known problem of
fictitious resonances when the wavenumber $k$ in the exterior problem
coincides with an eigenvalue of the interior problem. In this
approach, the Helmholtz equation is combined with its normal
derivative to yield a formulation which is numerically valid for all
real wavenumbers, at the expense of requiring the evaluation of
hypersingular integrals of the form
\begin{align*}
  \frac{\partial^{2}}{\partial n^{2}}
  \iint_{A_{e}} f(\xi,\eta)
  G(\mathbf{x},\mathbf{x}_{1}(\xi,\eta))\,\D \xi\,\D \eta.
\end{align*}
In order to meet continuity requirements, the collocation points in a
hypersingular method must lie strictly within elements, though
discontinuous elements offer a way around
this~\cite{marburg-schneider03}, and so we give a result for the
zero-order (constant) element only:
\begin{align}
  \label{equ:analysis:d2k:0}
  \frac{\partial^{2}K_{q,0}}{\partial z^{2}}
  =
  \frac{(kS)^{q}}{S}
  \biggr[
  \alpha
    &\int_{-\phi}^{\Theta-\phi}
    \left(
      \frac{\Delta}{\cos\theta} - \alpha
    \right)^{q+1}
    \frac{\cos^{3}\theta}{\Delta^{3}}
    \,\D\theta \\
    + (q+1)
    &\int_{-\phi}^{\Theta-\phi}
    \left(
      \frac{\Delta}{\cos\theta} - \alpha
    \right)^{q+1}
    \frac{\cos^{2}\theta}{\Delta^{2}}
    \,\D\theta
  \biggr].\nonumber
\end{align}

\subsection{Approximation of exponentials}
\label{sec:analysis:exp}

In order to efficiently evaluate the formulae of
Section~\ref{sec:analysis:basic}, we require a means of selecting the
polynomial approximation to the
exponential,~(\ref{equ:exp:approx}). The most obvious choice is to
truncate the Taylor series for $\E^{x}$ at some point where the
estimated remainder is smaller than a specified tolerance
$\epsilon$. For reasons of efficiency, however, we adopt an
``economized'' series which replaces the truncated Taylor series with
a polynomial approximation with minimum deviation and a minimized
error over the range where the polynomial is used. Given that the
integral terms are evaluated using recursion relations, by reducing
the number of terms, we also reduce the chance of numerical error
accumulating in moving from term to term.

The economization algorithm is that given by
Acton~\cite[p291--296]{acton90} and is used to generate a set of
polynomial approximations of $\sin x$ and $\cos x$ over a range $0\leq
x<\Delta x$, to a tolerance $\epsilon$ where $\epsilon$ is the maximum
difference between $\exp[\J x]$ and the polynomial approximation over the
range $0\leq x<\Delta x$. For the calculations of this paper, $\Delta
x=\pi/16,\pi/8,\pi/4,\pi/2$, and $\epsilon=10^{-n}$,
$n=3,6,9,12,15$. In the implementation, a polynomial approximation is
chosen which has the required maximum error less than $\epsilon$ and
$k\ell<\Delta x$. In the case of $\Delta x=\pi/2$, for example, this
gives a reduction in the number of terms required from fifteen for the
truncated Taylor series to eight for the economized polynomial when
$\epsilon=10^{-9}$.

\subsection{Summary of method}
\label{sec:analysis:summary}

The quadrature method of the previous sections can be summarized as
follows, for a triangle
$(\mathbf{x}_{1},\mathbf{x}_{2},\mathbf{x}_{3})$ which has been
rotated into the plane $z=0$ and field point $\mathbf{x}=(0,0,z)$,
Figure~\ref{fig:analysis:basic}:
\begin{enumerate}
\item determine the closest and furthest points on the triangle
  boundary and their radial distances $r_{\min}$
  (Figure~\ref{fig:analysis:basic}) and $r_{\max}$;
\item compute the required order of quadrature $Q$ for polynomial
  approximation of $1/R$, Section~\ref{sec:analysis:switch};
\item if $Q$ falls below the set limit:
  \begin{enumerate}
  \item  evaluate the integrals
  numerically and terminate;
  \end{enumerate}
  otherwise
  \begin{enumerate}
  \item decompose the triangle into up to three sub-triangles centered at
    the origin;
  \item for each sub-triangle compute the contribution using the
    formulae of Section~\ref{sec:analysis:basic} and accumulate,
    taking account of sub-triangle orientation.
  \end{enumerate}
\end{enumerate}

\section{Numerical testing}
\label{sec:testing}

\begin{figure}
  \centering
  \includegraphics{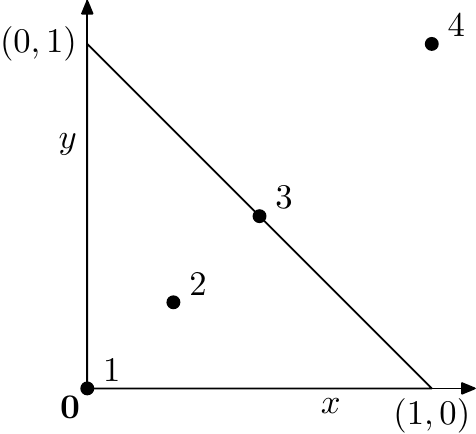}
  \caption{Sample triangle and field points}
  \label{fig:testing:tri}
\end{figure}

\begin{figure}
  \centering
  \begin{tabular}{l}
  \includegraphics{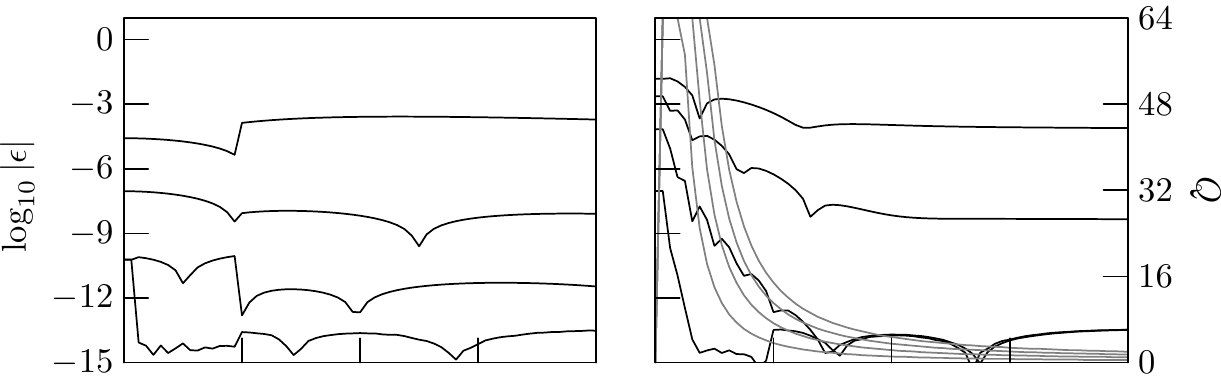} \\
  \includegraphics{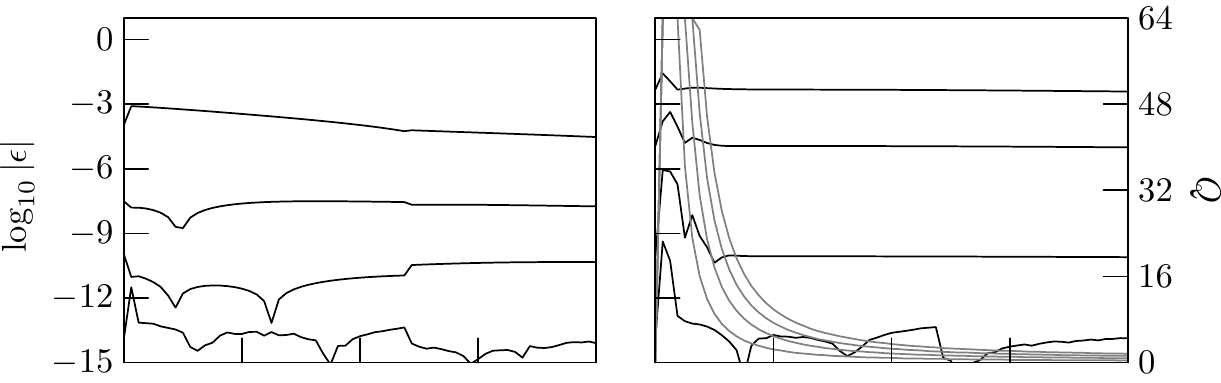} \\
  \includegraphics{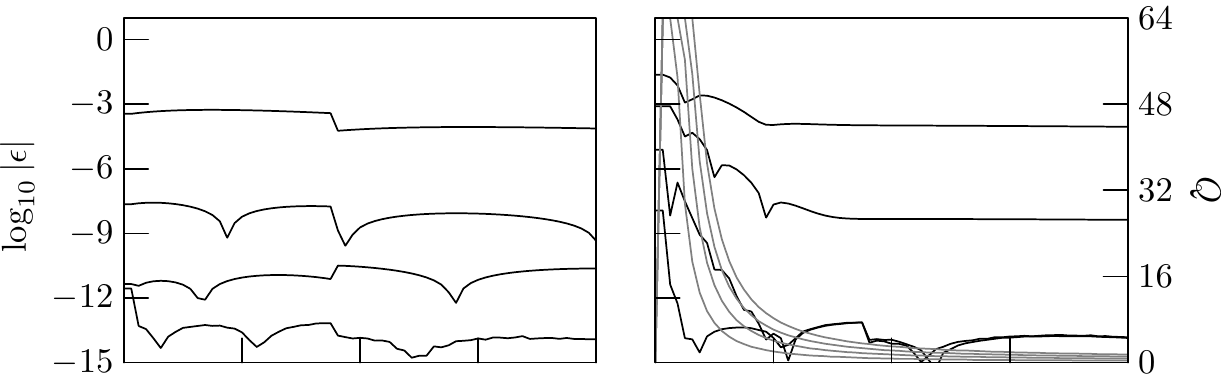} \\
  \includegraphics{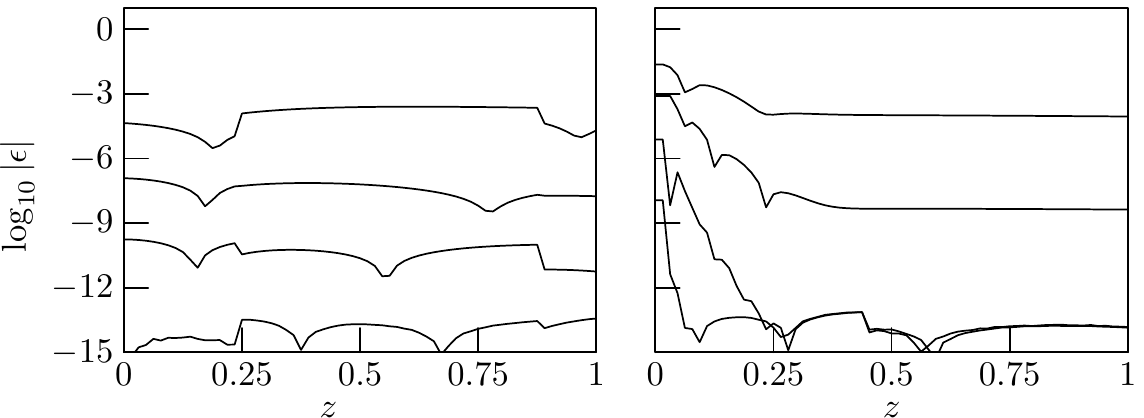}    
  \end{tabular}
  \caption{Error for $I_{0}(k,z)$ on element of
    Figure~\ref{fig:testing:tri} at points~1--4 (top to bottom);
    left-hand column: $|I_{0}^{(a)}-I_{0}^{(c)}|$ for
    $\epsilon^{(a)}=10^{-3,-6,-9,-12}$ and $50\times50$ point polar
    quadrature; right-hand column $|I_{0}^{(a)}-I_{0}^{(c)}|$ for
    $\epsilon^{(a)}=10^{-12}$ and $4\times4$, $8\times8$,
    $16\times16$, $32\times32$ polar quadrature; gray curves: $Q$ for
    $E_{Q}=10^{-3,-6,-9,-12}$}
  \label{fig:testing:p}
\end{figure}

\begin{figure}
  \centering
  \begin{tabular}{l}
  \includegraphics{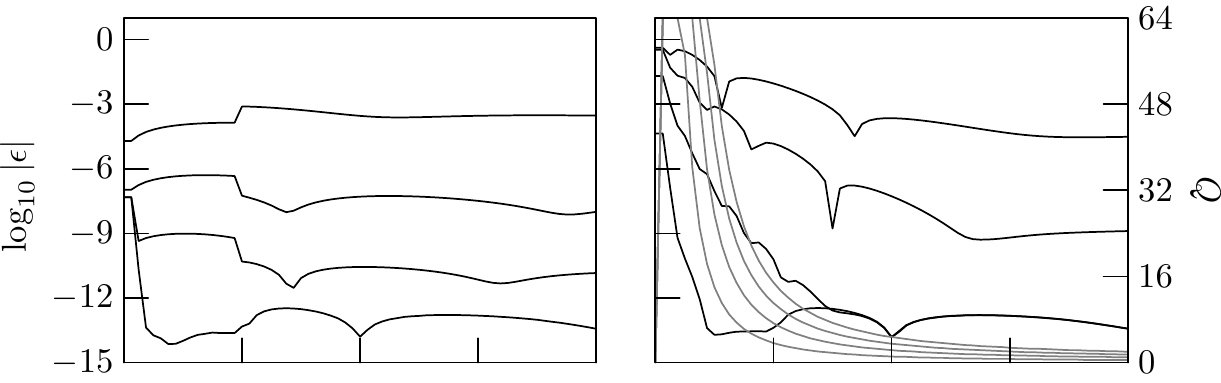} \\
  \includegraphics{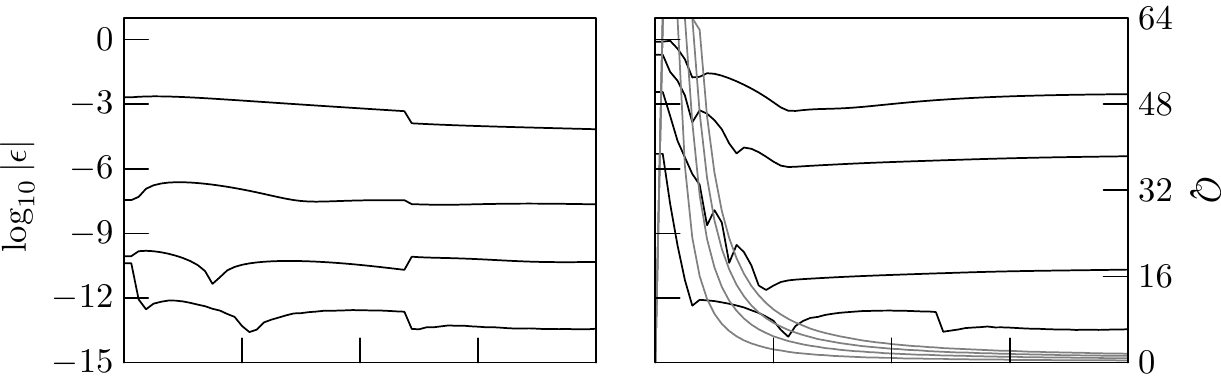} \\
  \includegraphics{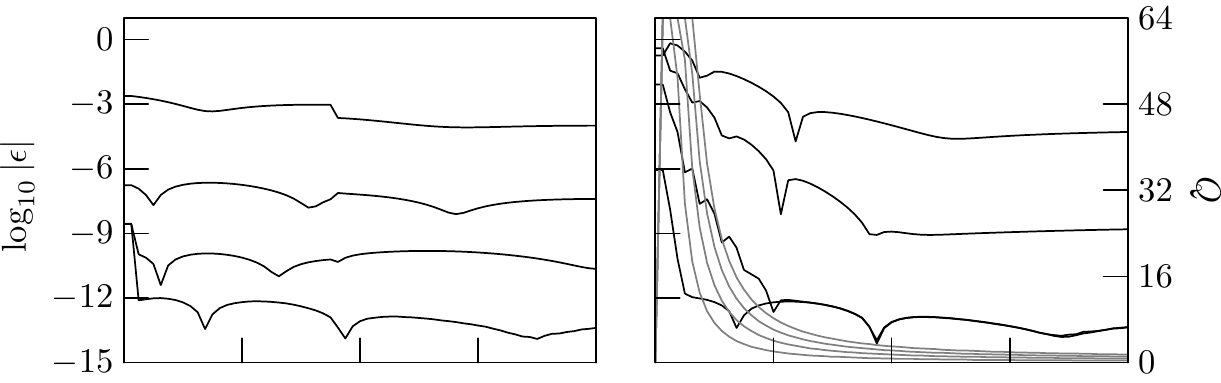} \\
  \includegraphics{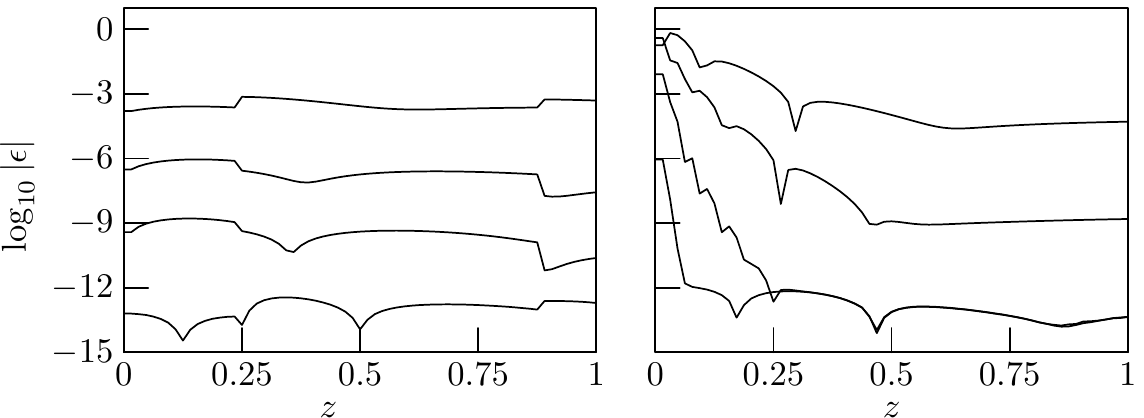}    
  \end{tabular}
  \caption{Error for $\partial I_{0}(k,z)/\partial n$ on element of
    Figure~\ref{fig:testing:tri} at points~1--4 (top to bottom),
    notation as in Figure~\ref{fig:testing:p}}
  \label{fig:testing:dp}
\end{figure}

As a numerical test of the performance of the method, we use the same
test case as Pourahmadian and
Mogilevskaya~\cite{pourahmadian-mogilevskaya15},
Figure~\ref{fig:testing:tri}. Four points are selected in the element
plane, as indicated, and we evaluate $I_{0}$ for $k=1$ as a function
of $z$, vertical displacement from the element; results for $I_{x}$
and $I_{y}$ are similar. As a reference for error estimation, we
follow Pourahmadian and Mogilevskaya and use a polar transformation
and a $50\times50$ point Gaussian quadrature, which is accurate to
eight significant figures~\cite{pourahmadian-mogilevskaya15}. The
reported error is
\begin{align}
  \epsilon(z) &= |I_{0}^{(a)}(z)-I_{0}^{(c)}(z)|,
\end{align}
with superscripts `a' and `c' denoting `analytical' and `computed'
values respectively. Error is evaluated by specifying the required
tolerance $\epsilon^{(a)}$ in the analytical method and computing the
resulting $\epsilon$. A second set of error calculations are presented
by fixing $\epsilon^{(a)}=10^{-12}$, varying the order of numerical
quadrature in the polar transformation, and computing the resulting
error $\epsilon^{(c)}$. Figure~\ref{fig:testing:p} gives $\epsilon$ as
a function of $z$ for varying $\epsilon^{(a)}$, and $\epsilon^{(c)}$
for varying order of Gaussian quadrature, plotted with $Q$ computed
for varying values of $E_{Q}$. Figure~\ref{fig:testing:p} shows error
data for the evaluation of $I_{0}$ and Figure~\ref{fig:testing:dp} for
the normal derivative $\partial I_{0}/\partial z$.

The left-hand column of Figure~\ref{fig:testing:p} shows the error
estimate for points~1,~2,~3, and~4 in Figure~\ref{fig:testing:tri}
which correspond respectively to field points whose projections lie on
a vertex, in the interior, on an edge, and outside the element. The
reference integral $I_{0}^{(c)}$ for points~1--3 is computed using the
$50\times50$ Gaussian quadrature after transformation to polar
coordinates, and that for point~4 using the~175 point symmetric
quadrature of Wandzura and Xiao~\cite{wandzura-xiao03}. Errors are
computed with a requested tolerance $\epsilon^{(a)}=10^{-3,-6,-9,-12}$
and the computed errors reflect both the accuracy of the analytical
method and the conformity to the requested tolerance. In applications,
there is reasonable confidence that the error will be approximately
equal to that requested, without excessive computation.

The right-hand column of Figure~\ref{fig:testing:p} presents data
relevant to the error estimate $E_{Q}$ and the accuracy of Gaussian
quadrature in this problem. The darker curves show an error estimate
computed as the difference between the analytical method with a
requested tolerance of $10^{-12}$ and numerical quadrature of varying
order. As expected the low-order methods, e.g. $4\times4$ points, give
a larger error and the high-order approach, $32\times32$ points, gives
accuracy comparable to the analytical technique, except for small
values of $z$. In each case, the breakdown of the polynomial
approximation for $1/R$ is apparent in the increase of error as
$z\to0$, most clearly for the $16\times16$ quadrature where the error
increases markedly from $z\approx0.3$. Of interest here is the value
of $E_{Q}$ as a criterion for selecting quadrature rules. The lighter
curves show the value of $Q$ found from~(\ref{equ:error:abs})
with varying values of $E_{Q}$. The curves do indeed predict quite
well the point at which the polynomial approximation to $1/R$ is no
longer accurate and the Gaussian quadrature begins to fail, confirming
the reliability of the measure as a criterion for the selection of
quadrature rules.

Figure~\ref{fig:testing:dp} gives similar results but for the
evaluation of the normal derivative of the layer potential, also
required in BEM calculations. The results are similar to those in
Figure~\ref{fig:testing:p} and the discussion of those data carries
over to here, but it is worth noting that though the error behavior of
the Gaussian quadratures is different from that in
Figure~\ref{fig:testing:p} (compare the results for point~3, for
example), the curves of $Q$ still function as a reliable criterion for
selecting a quadrature method.

\section{Conclusions}
\label{sec:conclusions}

An analytical method for the evaluation of potential integrals in
boundary element codes for the Helmholtz equation has been presented
and tested. An error estimator for purely numerical quadrature has
been derived and used to establish a criterion for quadrature method
selection. The quadrature method has been tested and found to be
accurate and reliable; the error criterion is a reliable technique for
quadrature selection. We believe that the quadrature method proposed
is a suitable plug-in replacement in BEM codes for the wave equation
where an \textit{a priori} error estimate for element integrals and an
economical integration are required.

\appendix

\section{Basic integrals}
\label{sec:basic:integrals}

The evaluation of the potential integrals requires a number of
elementary integrals which can be computed using results from standard
tables combined with recursions. This appendix contains the results
required for the evaluation of the trigonometric integrals of the main
paper, written in terms of the parameter $\alpha$, $0\leq\alpha<1$,
and $\Delta^{2}=1-\alpha^{2}\sin^{2}\theta$. The results are given as
the indefinite integral, with a separate result where necessary for
the in-plane case $\alpha=0$. 

The first basic term is
\begin{align}
  \int
  \left(
    \frac{\Delta}{\cos\theta}
    -
    \alpha
  \right)^{q}
  \left(
    \frac{\Delta}{\cos\theta}  
  \right)^{-s}
  \,\D\theta
  &=
  \sum_{u=0}^{q}
  \binom{q}{u}
  \left(-\alpha\right)^{q}
  \int
  \left(
    \frac{\Delta}{\cos\theta}
  \right)^{q-s}
  \,\D\theta,
\end{align}
where $s=0,1,2,3$. The terms in the summation are pseudo-elliptic
integrals which can be evaluated using elementary functions and
recursion relations~\cite[2.58]{gradshteyn-ryzhik80}.

Using the transformation $u=\tan\theta$ and noting that
$\Delta/\cos\theta=(1+{\alpha'}^{2}\tan^{2}\theta)^{1/2}$
\begin{align}
  \label{equ:basic:one}
  \int 
  \left(\frac{\Delta}{\cos\theta}\right)^{n}
  \,\D\theta
  &=
  \alpha^{2}
  \int 
  \left(\frac{\Delta}{\cos\theta}\right)^{n-2}
  \,\D\theta
  + {\alpha'}^{2}
  \int 
  \left(1+{\alpha'}^{2}u^{2}\right)^{(n-2)/2}
  \,\D u,
\end{align}
with $\alpha' = (1-\alpha^{2})^{1/2}$.

The integral term can be evaluated using the recursion
\begin{align}
  \int 
  \left(1+{\alpha'}^{2}u^{2}\right)^{(n+2)/2}
  \,\D u 
  &=
  \frac{u}{n+3}\left(1+{\alpha'}^{2}u^{2}\right)^{(n+2)/2}
  +
  \frac{n+2}{n+3}
  \int \left(1+{\alpha'}^{2}u^{2}\right)^{n/2}\,\D u,
\end{align}
seeding the recursion with
\begin{align}
  \int \left(1+{\alpha'}^{2}u^{2}\right)^{-2/2}\,\D u
  &=
  \frac{\tan^{-1}(\alpha' u)}{\alpha'},\\
  \int \left(1+{\alpha'}^{2}u^{2}\right)^{-1/2}\,\D u
  &=
  \frac{1}{\alpha'}\log \left[
    \left(1+{\alpha'}^{2}u^{2}\right)^{1/2} + \alpha'u
  \right],
\end{align}
and using
\begin{align}
  \int 
  \left(\frac{\Delta}{\cos\theta}\right)^{-3}
  \,\D\theta
  &=
  -\frac{{\alpha'}^{2}}{\alpha^{2}}\frac{\sin\theta}{\Delta}
  + \frac{\sin^{-1}(\alpha\sin\theta)}{\alpha},\\
  \int 
  \left(\frac{\Delta}{\cos\theta}\right)^{-2}
  \,\D\theta
  &=
  \frac{\theta}{\alpha^{2}}
  -\frac{\alpha'}{\alpha^{2}}
  \tan^{-1}(\alpha'\tan\theta),
  \\
  \int 
  \left(\frac{\Delta}{\cos\theta}\right)^{-1}
  \,\D\theta
  &=
  \frac{\sin^{-1}(\alpha\sin\theta)}{\alpha}.
\end{align}

For $\alpha=0$,
\begin{align}
  \int 
  \left(\frac{\Delta}{\cos\theta}\right)^{-3}
  \,\D\theta
  &=
  \sin\theta - \frac{\sin^{3}\theta}{3},\\
  \int 
  \left(\frac{\Delta}{\cos\theta}\right)^{-2}
  \,\D\theta
  &=
  \frac{\sin\theta\cos\theta}{2} + \frac{\theta}{2}
  \\
  \int 
  \left(\frac{\Delta}{\cos\theta}\right)^{-1}
  \,\D\theta
  &=
  \sin\theta.
\end{align}

A second, similar, integral is
\begin{align}
  \label{equ:basic:tan}
  \int 
  \left(\frac{\Delta}{\cos\theta}\right)^{n}
  \frac{\sin\theta}{\cos\theta}
  \,\D\theta
  &=
  \alpha^{2}
  \int 
  \left(\frac{\Delta}{\cos\theta}\right)^{n-2}
  \frac{\sin\theta}{\cos\theta}
  \,\D\theta
  +
  \frac{1}{n}
  \left(
    \frac{\Delta}{\cos\theta}
  \right)^{n},
\end{align}
which can be seeded with~\cite[2.584]{gradshteyn-ryzhik80}
\begin{align*}
  \int 
  \left(\frac{\Delta}{\cos\theta}\right)^{-3}
  \frac{\sin\theta}{\cos\theta}
  \,\D\theta
  &= 
  \frac{\cos\theta}{\alpha^{2}\Delta}
  -
  \frac{1}{\alpha^{3}}
  \log\left(\alpha\cos\theta + \Delta\right),\\
  \int 
  \left(\frac{\Delta}{\cos\theta}\right)^{-2}
  \frac{\sin\theta}{\cos\theta}
  \,\D\theta
  &=
  -\frac{1}{\alpha^{2}}\log\Delta,\\
  \int 
  \left(\frac{\Delta}{\cos\theta}\right)^{-1}
  \frac{\sin\theta}{\cos\theta}
  \,\D\theta
  &= 
  -\frac{1}{\alpha}\log(\alpha\cos\theta + \Delta),\\
  \int 
  \left(\frac{\Delta}{\cos\theta}\right)^{0}
  \frac{\sin\theta}{\cos\theta}
  \,\D\theta
  &=
  -\log\cos\theta.
\end{align*}

For $\alpha=0$,
\begin{align*}
  \int 
  \left(\frac{\Delta}{\cos\theta}\right)^{-3}
  \frac{\sin\theta}{\cos\theta}
  \,\D\theta
  &= 
  -\frac{\cos^{3}\theta}{3},\\
  \int 
  \left(\frac{\Delta}{\cos\theta}\right)^{-2}
  \frac{\sin\theta}{\cos\theta}
  \,\D\theta
  &=
  \frac{\sin^{2}\theta}{2},\\
  \int 
  \left(\frac{\Delta}{\cos\theta}\right)^{-1}
  \frac{\sin\theta}{\cos\theta}
  \,\D\theta
  &= 
  -\cos\theta.
\end{align*}

In an implementation of the method of this paper, when the required
geometric parameters have been calculated for the reference triangle,
and the appropriate expansion for $\exp[\J k x]$ has been selected,
the first step is to compute the required elementary integrals
(\ref{equ:basic:one}) and (\ref{equ:basic:tan}) using the initial
values and the recursion relations. The computed terms can then be
used in the summations of (\ref{sec:analysis:basic}) to evaluate the
potential integrals.

For convenience, we define
\begin{align}
  \label{equ:integrals:log:cos}
  L_{c} &=
  \int 
  \cos\theta
  \log\frac{\Delta-\alpha'}{\Delta+\alpha'}\,\D\theta
  =
  \sin\theta
  \log\frac{\Delta-\alpha'}{\Delta+\alpha'}
  +\log\frac{\Delta+\alpha'\sin\theta}{\Delta-\alpha'\sin\theta}
  - 2\frac{\alpha'}{\alpha}
  \sin^{-1}(\alpha\sin\theta),\\
  \label{equ:integrals:log:sin}
  L_{s} &=
  \int 
  \sin\theta
  \log\frac{\Delta-\alpha'}{\Delta+\alpha'}\,\D\theta
  =
  -\cos\theta
  \log\frac{\Delta-\alpha'}{\Delta+\alpha'}
  + 2\frac{\alpha'}{\alpha}
  \log(\alpha\cos\theta+\Delta),
\end{align}
which are readily evaluated using integration by parts. Then,
\begin{align}
  I_{q,c} &= \int J_{q} \cos\theta\,\D\theta,\\
  I_{q+1,c} &= \frac{s(kS)^{q+1}}{2(q+2)}
  \int \left(\frac{\Delta}{\cos\theta} - \alpha\right)^{q+1}\,\D\theta
  -k|z|\frac{2q+3}{q+2}I_{q,c},\\
  I_{0,c} &= \frac{s}{2}\theta + \frac{|z|}{4}L_{c},\\
  \frac{\partial I_{q+1,c}}{\partial z}
  &=
  \mp \frac{s}{2S}(kS)^{q+1} \frac{q+1}{q+2}
  \int 
  \left(\frac{\Delta}{\cos\theta} - \alpha\right)^{q+1}
  \frac{\cos\theta}{\Delta}
  \,\D\theta\nonumber\\
  &\mp
  \frac{2q+3}{q+2}
  k I_{q,c}
  -k|z|
  \frac{2q+3}{q+2}
  \frac{\partial I_{q,c}}{\partial z},\\
  \frac{\partial I_{0,c}}{\partial z}
  &=
  \pm L_{c} \pm \frac{s}{2S}
  \int
  \frac{\cos\theta}{\Delta}
  \,\D\theta,
\end{align}
and similarly
\begin{align}
  I_{q+1,s} &= \frac{s(kS)^{q+1}}{2(q+2)}
  \int 
  \left(\frac{\Delta}{\cos\theta} - \alpha\right)^{q+1}
  \frac{\sin\theta}{\cos\theta}
  \,\D\theta
  -k|z|\frac{2q+3}{q+2}I_{q,s},\\
  I_{0,s} &= -\frac{s}{2}\log\cos\theta + \frac{|z|}{4}L_{s},\\
  \frac{\partial I_{q+1,s}}{\partial z}
  &=
  \mp \frac{s}{2S}(kS)^{q+1} \frac{q+1}{q+2}
  \int 
  \left(\frac{\Delta}{\cos\theta} - \alpha\right)^{q+1}
  \frac{\sin\theta}{\Delta}
  \,\D\theta\nonumber\\
  &\mp
  \frac{2q+3}{q+2}
  k I_{q,s}
  -k|z|
  \frac{2q+3}{q+2}
  \frac{\partial I_{q,s}}{\partial z},\\
  \frac{\partial I_{0,s}}{\partial z}
  &=
  \pm L_{s} \pm \frac{s}{2S}
  \int
  \frac{\sin\theta}{\Delta}
  \,\D\theta
\end{align}

\bibliographystyle{siamplain}

\end{document}